\documentclass{conm-p-l}
\usepackage{amssymb,amsmath,amsthm}

\numberwithin{equation}{section}

\usepackage{enumerate}

\newtheorem{theorem}{Theorem}[section]

\newtheorem{corollary}[theorem]{Corollary}

\theoremstyle{definition}

\newtheorem{definition}[theorem]{Definition}
\newtheorem{notation}[theorem]{Notation}
\newtheorem{example}[theorem]{Example}

\newtheorem{remark}[theorem]{Remark}

\newtheorem{conjecture}[theorem]{Conjecture}

\newcommand{\C}{{\mathbb{C}}}
\newcommand{\N}{{\mathbb{N}}}
\newcommand{\R}{{\mathbb{R}}}

\newcommand{\Z}{{\mathbb{Z}}}

\newcommand{\cor}{{\bf{k}}}
\newcommand{\cora}{{\BBK}}

\newcommand{\coro}{{{\bf{k}}_0}}
\newcommand{\hcor}{{\bf{\widehat k}}}

\newcommand{\rmw}{{\rm w}}
\newcommand{\SSi}{\mathrm{SS}}

\def\BBK{\mathbb{K}}

\def\phi{{\varphi}}
\def\epsilon{\varepsilon}

\def\sha{\mathcal{A}}

\def\shd{\mathcal{D}}
\def\she{\mathcal{E}}

\def\shi{\mathcal{I}}

\def\shk{\mathcal{K}}
\def\shl{\mathcal{L}}
\def\shm{\mathcal{M}}
\def\shn{\mathcal{N}}
\def\sho{\mathcal{O}}

\def\sht{\mathcal{T}}
\def\shu{\mathcal{U}}

\def\shw{\mathcal{W}}

\newcommand{\OO}{\sho}

\newcommand{\atw}[1]{{}_{#1}}

\newcommand{\halfform}{{\sqrt v}}

\newcommand{\E}[1][]{\mathcal{E}_{#1}}
\newcommand{\HE}[1][]{\widehat{\mathcal{E}}_{#1}}

\newcommand{\Ev}[1][]{\mathcal{E}_{#1}^\halfform}

\newcommand{\W}[1][]{\mathcal{W}_{#1}}

\newcommand{\Wva}[1][]{\mathcal{W}_{#1}^{\halfform,+}}

\newcommand{\CD}[1][]{\mathcal{C}_{#1}}

\newcommand{\HW}[1][]{\widehat{\mathcal{W}}_{#1}}
\newcommand{\SW}[1][]{\mathcal{W}^s_{#1}}

\newcommand{\TW}[1][]{\mathcal{W}^t_{#1}}

\newcommand{\Wv}[1][]{\mathcal{W}_{#1}^\halfform}

\newcommand{\gr}{\mathop{\mathrm{gr}}}

\newcommand{\Lie}[1][]{\operatorname{\mathsf{L}}\def\temp{#1}
\ifx\temp\empty\else^{(#1)}\fi}
\newcommand{\rmad}{\mathrm{Ad}}

\newcommand{\stx}{{\mathfrak{X}}}
\newcommand{\sty}{{\mathfrak{Y}}}

\newcommand{\stkHom}[1][]{\mathfrak{Hom}_{\raise1.5ex\hbox to.1em{}#1}}

\newcommand{\rmpt}{{\rm pt}}
\newcommand{\rmptt}{{\{\rm pt\}}}
\newcommand{\Eu}{{\rm Eu}}
\newcommand{\Ch}{{\rm Ch}}
\newcommand{\Td}{{\rm Td}}

\renewcommand{\to}[1][]{\xrightarrow[]{#1}}

\newcommand{\isoto}[1][]{\xrightarrow[#1]%
{{\raisebox{-.6ex}[0ex][-.6ex]{$\mspace{1mu}\sim\mspace{2mu}$}}}}
\newcommand{\isofrom}[1][]{\xleftarrow[#1]{\sim}}

\DeclareMathOperator{\id}{id}

\newcommand{\Hom}[1][]{\mathrm{Hom}_{\raise1.5ex\hbox to.1em{}#1}}
\newcommand{\RHom}[1][]{\mathrm{RHom}_{\raise1.5ex\hbox to.1em{}#1}}
\newcommand{\Ext}[2][]{\mathrm{Ext}_{\raise1.5ex\hbox to.1em{}#1}^{#2}}
\renewcommand{\hom}[1][]{{\mathcal{H}om}_{\raise1.5ex\hbox to.1em{}#1}}
\newcommand{\rhom}[1][]{{R\mathcal{H}om}_{\raise1.5ex\hbox to.1em{}#1}}
\newcommand{\ext}[2][]{{\mathcal{E}xt}_{\raise1.5ex\hbox to.1em{}#1}^{#2}}

\newcommand{\Tens}[1][]{\mathbin{\otimes_{\raise1.5ex\hbox to-.1em{}#1}}}
\newcommand{\LTens}[1][]{\mathbin{\otimes_{\raise1.5ex\hbox to-.1em{}#1}^{L}}}
\newcommand{\Tor}[2][]{\mathrm{Tor}^{\raise1.5ex\hbox to.1em{}#1}_{#2}}
\newcommand{\tens}[1][]{\mathbin{\otimes_{\raise1.5ex\hbox to-.1em{}#1}}}
\newcommand{\ltens}[1][]{\mathbin{\otimes_{\raise1.5ex\hbox to-.1em{}#1}^{L}}}
\newcommand{\detens}{\underline{\etens}}
\newcommand{\tor}[2][]{{\mathcal{T}or}^{\raise1.5ex\hbox to.1em{}#1}_{#2}}
\newcommand{\etens}{\mathbin{\boxtimes}}

\newcommand{\lltens}[1][]{{\mathop{\tens}\limits^{\rm L}}_{#1}}

\newcommand{\Endo}[1][]{\mathrm{End}_{\raise1.5ex\hbox to.1em{}#1}}

\newcommand{\Aut}[1][]{\mathrm{Aut}_{\raise1.5ex\hbox to.1em{}#1}}

\newcommand{\RD}{{\rm D}}

\newcommand{\Rb}{{\rm b}}
\newcommand{\coh}{{\rm coh}}

\newcommand{\gd}{{\rm gd}}
\newcommand{\gdc}{{\rm gd,c}}

\newcommand{\oim}[1]{{#1}_*}
\newcommand{\eim}[1]{{#1}_!}

\newcommand{\reim}[1]{{R#1}_!}
\newcommand{\opb}[1]{#1^{-1}}

\DeclareMathOperator{\supp}{supp}
\DeclareMathOperator{\ori}{or}

\def\rop{{\rm op}}
\def\tot{{\rm tot}}

\newcommand{\res}{{\rm res}}
\newcommand{\hol}{{\rm hol}}
\newcommand{\rh}{{\rm rh}}
\newcommand{\Cc}{{\C\rm c}}
\newcommand{\Rc}{{\R\rm c}}

\newcommand{\md[1]}{{\rm Mod}({#1})}

\newcommand{\Mods}{\operatorname{{\mathfrak{Mod}}}}

\newcommand{\eqdot}{\mathbin{:=}}

\newcommand{\cl}{\colon}
\newcommand{\scbul}{\,\raise.4ex\hbox{$\scriptscriptstyle\bullet$}\,}

\newcommand{\ba}{\begin{array}}
\newcommand{\ea}{\end{array}}

\newcommand{\bnum}{\begin{enumerate}[{\rm(i)}]}
\newcommand{\enum}{\end{enumerate}}
\newcommand{\banum}{\begin{enumerate}[{\rm(a)}]}
\newcommand{\eanum}{\end{enumerate}}

\newcommand{\lp}{{\rm(}}
\newcommand{\rp}{{\rm)}}

\newcommand{\eq}{\begin{eqnarray}}
\newcommand{\eneq}{\end{eqnarray}}
\newcommand{\eqn}{\begin{eqnarray*}}
\newcommand{\eneqn}{\end{eqnarray*}}

\begin{document}

\title{Deformation quantization modules on complex symplectic manifolds}
\author{Pierre Schapira} 
\address{
Institut de Math{\'e}matiques,
Universit{\'e} Pierre et Marie Curie,
175, rue du Chevaleret, 75013 Paris, France
}
\email {schapira@math.jussieu.fr}
\urladdr{http://www.math.jussieu.fr/$\sim$schapira}




\maketitle

\begin{abstract}
We study modules over the algebroid stack $\W[\stx]$ of
deformation quantization on a complex symplectic manifold $\stx$ and
recall some results: construction of an algebra for $\star$-products,
existence of (twisted) simple modules along smooth
Lagrangian submanifolds, perversity of the complex of solutions for
regular holonomic  $\W[\stx]$-modules, finiteness and duality for
the composition of ``good''  kernels. 
As a corollary, we get that the derived category
of good $\W[\stx]$-modules with compact support is a Calabi-Yau
category. We also give a conjectural 
Riemann-Roch type formula in this framework.
\end{abstract}

\section*{Introduction}
Let $X$ be  a complex manifold, $T^*X$ its cotangent bundle.
The conic sheaf of $\C$-algebras  $\E[T^*X]$ 
of microdifferential operators on 
$T^*X$  has been constructed
functorially by Sato-Kashiwara-Kawai in \cite{S-K-K}. 
This algebra is associated with the homogeneous symplectic
structure and is also naturally defined on the projective cotangent
bundle $P^*X$. 

Another (no more conic) algebra on $T^*X$,
denoted here by  $\W[T^*X]$ and defined over
a subfield $\cor$ of $\C[[\opb{\tau},\tau]$
has been constructed in \cite{P-S} (see \cite{Bo} for related
constructions). 
Its formal version has been
considered by  many authors after \cite{BFFLS} and extended to Poisson
manifolds in \cite{Ko2}. 

In general, neither the algebras  $\E[P^*X]$ glue on a complex 
contact manifold, nor the algebras $\W[T^*X]$  glue on a 
complex symplectic manifold, although the
 categories of modules on these non existing algebras make sense. 
Indeed, one has to replace the notion of a sheaf of algebras by that
 of an algebroid stack, similarly as one replaces the notion of a  sheaf by
 that of a stack. These constructions are  performed 
in \cite{Ka1}, \cite{Ko}, \cite{P-S} (see also \cite{B-G-N-T} for
recent developments and \cite{Bz-K,VdB,Ye} for an algebraic approach).

Here, we start by briefly recalling the constructions of the sheaves $\E[T^*X]$
and $\W[T^*X]$ as well as a new sheaf of algebras on 
$T^*X$ containing $\W[T^*X]$,
invariant by quantized symplectic transformations, in which the
$\star$-exponential is well defined (see \cite{Di-S}).
Then we consider a complex symplectic manifold $\stx$,
introduce the algebroid stack $\W[\stx]$ of deformation
quantization on $\stx$ and discuss some recent
 results on $\W[\stx]$-modules:
\begin{itemize}
\item
If $\Lambda$ is a smooth Lagrangian submanifold of $\stx$, there
exist twisted simple $\W[\stx]$-modules along $\Lambda$, 
the twist being associated
with a square root of the line bundle  $\Omega_\Lambda$ (see \cite{Da-S}).
\item
Let $\shl_0$ and $\shl_1$ be two regular
holonomic modules supported by smooth Lagrangian submanifolds
$\Lambda_0$ and $\Lambda_1$. 
Then the complex $\rhom[{\W[\stx]}](\shl_0,\shl_1)$ is a perverse sheaf over the
field $\cor$  (see \cite{K-S3}). 
\item
Let $\stx_i$ ($i=1,2,3$) be complex symplectic manifolds 
 and denote by $\stx_i^a$  the symplectic manifold deduced from $\stx_i$ 
by taking the opposite symplectic form. Let 
$\shk_i$ be a good $\W[\stx_{i+1}\times\stx_i^a]$-module ($i=1,2$) 
(good means coherent and endowed with a good filtration on
each compact subset of $\stx$) and assume that 
a properness condition is satisfied by the supports of these
modules. Then their composition $\shk_2\circ\shk_1$ is  a good 
$\W[\stx_{3}\times\stx_1^a]$-module. Moreover, composition of kernels
commutes with duality (see \cite{Sc-S2}). 
As a particular case, we obtain that the triangulated category 
consisting of good $\W[\stx]$-modules  with compact supports 
is   Ext-finite over the field $\cor$ and admits a Serre
functor, namely the shift by $d_\stx\eqdot\dim_\C\stx$. 
\item 
The Hochschild homology of the algebroid stack $\W[\stx]$ is
concentrated in degree $-\dim_\C\stx$ and is isomorphic to $\cor_\stx$.
This allows us to construct the
Euler class $\Eu(\shm)\in H^{d_\stx}_{\supp\shm}(\stx;\cor_\stx)$ of
a coherent $\W[\stx]$-module. We conjecture that in the situation 
above, $\Eu(\shk_2\circ\shk_1)$ is the relative integral of the cup
products $\Eu(\shk_1)\cup\Eu(\shk_2)$ (see \cite{Sc-S2}). 
\end{itemize}
This paper summarizes various joint works with 
A.~D'Agnolo \cite{Da-S}, G.~Dito \cite{Di-S}, M.~Kashiwara \cite{K-S3}, 
P.~Polesello \cite{P-S} and J-P.~Schneiders \cite{Sc-S2}.

\section{Microdifferential operators on cotangent bundles}

Let $X$ be a complex manifold, $\pi\cl T^*X\to X$ its cotangent bundle.

\subsection*{The ring $\she_{T^*X}$.}

The manifold $T^*X$ is a complex {\em homogeneous} symplectic 
manifold, {\em i.e.,}  
$T^*X$ is endowed with a canonical $1$-form $\alpha_X$ such that   
$d\alpha_X$ is symplectic.  
On $T^*X$ there exists a conic ({\em i.e.,} constant on the orbits of the action
of $\C^\times$)  sheaf
 $\E[T^*X]$ constructed functorially by Sato-Kashiwara-Kawai \cite{S-K-K} 
(see also \cite{Ka2,S} for an exposition)
which plays the role of a noncommutative 
localization of the ring $\shd_X$ of differential operators. 
The sheaf $\E[T^*X]$ enjoys the following properties:
\begin{itemize}
\item
$\E[T^*X]$ is a filtered sheaf of central $\C$-algebras and 
\eqn
&&\gr\E[T^*X]\simeq \bigoplus_{j\in\Z}\sho_{T^*X}(j).
\eneqn
($\sho_{T^*X}(j)$  is the subsheaf of $\sho_{T^*X}$ consisting of
homogeneous functions of degree $j$ in the fibers of $\pi$.)
\item
There is a flat monomorphism of filtered rings
$\opb{\pi}\shd_X\hookrightarrow\E[T^*X]$. 
\item
Denote by $\Omega_X^{\frac{1}{2}}$ 
the twisted sheaf of holomorphic  half-forms of maximal degree on $X$ 
(the notion of twisted sheaves will be recalled below). 
One defines the sheaf of algebras:
\eq\label{def:Ev}
&&\Ev[T^*X]\eqdot
\opb{\pi}\Omega_X^{\frac{1}{2}}\tens[\opb{\pi}\sho_X]\E[T^*X]
\tens[\opb{\pi}\sho_X]\opb{\pi}\Omega_X^{-\frac{1}{2}}.
\eneq
(Note that $\Ev[T^*X]$ is a sheaf although $\Omega_X^{\frac{1}{2}}$ 
is not a sheaf but a twisted sheaf.) 
Denote by $a\cl T^*X\to T^*X$ the antipodal map, 
$(x;\xi)\mapsto (x;-\xi)$. 
There exists a $\C$-linear anti-isomorphism of sheaves of 
algebras\footnote{An anti-isomorphism of algebras $A\isoto B$ is an
  isomorphism of algebras $A\isoto B^\rop$ where 
$B^\rop$ is the opposite algebra.}
$\oim{a}\Ev[T^*X]\isoto \Ev[T^*X]$ called the transposition and
denoted  $P\mapsto {}^tP$. 
\item
Consider a $\C^\times$-homogeneous symplectic isomorphism
$\phi\cl T^*X\supset U\isoto V\subset T^*Y$. 
Then $\phi$ can be {\em locally} quantized as an 
isomorphism of filtered sheaf of rings commuting with the transposition
\eqn
&& \Phi\cl \oim{\phi}\Ev[T^*X]\isoto \Ev[T^*Y].
\eneqn
Denote by $V^a$ the image of $V$ by the antipodal map $a$ on $T^*Y$
and by $\Lambda_\phi\subset U\times V^a$ the image of the graph of
$\phi$. This is a Lagrangian submanifold of $U\times V^a$.  
Locally, we may assume that $\Omega_X$ is trivial and 
there exists an ideal $\shi_\phi$
of $\E[T^*(X\times Y)]$ whose  associated graded
ideal is reduced and coincides with the defining ideal of $\Lambda_\phi$.
Then, for each 
$\E[T^*X]\ni P$ there exists a unique $Q\in\E[T^*Y]$ such that
$P-Q\in\shi_\phi$. 
The correspondence $P\mapsto Q$ is an anti-isomorphism of $\C$-algebras
$\oim{\phi}\E[T^*X]\isoto \oim{a}\E[T^*Y]$. One gets the isomorphism
$\Phi$ by composing with the transposition in $\E[T^*Y]$. 
One shall be aware that 
this isomorphism exists {\em only locally and is not unique} in general.
\end{itemize}
Moreover, when $X$ is affine ({\em i.e.,} 
open in some $n$-dimensional complex vector
space), $\E[T^*X]$ satisfies:
\begin{itemize}
\item
 any section $P\in\E[T^*X](U)$ on an open subset $U\subset T^*X$ admits a 
total symbol
\eq\label{eq:totsymbe}
&&\sigma_{\tot}(P)(x;\xi) = \sum_{-\infty<j\leq m}p_j(x;\xi), 
\,m\in\Z\, \quad p_j\in\sho_{T^*X}(j)(U),
\eneq
with the condition:
\eq\label{eq:defE}
&&\left\{  \parbox{300 pt}{
for any compact subset $K$ of $U$ there exists a positive 
constant $C_K$ such that 
$\sup\limits_{K}\vert p_{j}\vert \leq C_K^{-j}(-j)!$ for all $j\leq 0$.
}\right. 
\eneq
\item
The total symbol of the product is given by the Leibniz rule:
\eq\label{eq:leibniz2}
\sigma_{\tot}(P\circ Q)
&=&\sum_{\alpha\in\N^n} \dfrac{1}{\alpha !} 
\partial^{\alpha}_\xi\sigma_{\tot}(P)\partial^{\alpha}_x\sigma_{\tot}(Q).
\eneq
\item
The total symbol of the transposition is given by 
\eq\label{eq:transposeE}
\sigma_{\tot}({}^tP)(x;\xi)
&=&\sum_{\alpha\in\N^n}\dfrac{(-1)^{\vert\alpha\vert}}{\alpha!}
\partial_\xi^\alpha\partial_x^\alpha\sigma_{\tot}(P)(x;-\xi).
\eneq
\end{itemize}

\subsection*{The field $\cor$.}

Let $\hcor\eqdot \C[[\opb{\tau},\tau]$ be the field of formal Laurent
series in $\opb{\tau}$.
We consider the filtered subfield $\cor$ of $\hcor$
consisting of series
$a = \sum_{-\infty< j\leq m}a_j\tau^j$ ($a_j\in\C$,  $m\in\Z$) satisfying:
\eq\label{eq:cor}
&&\parbox{300 pt}{
there exists 
$C>0$ such that $\vert a_{j}\vert \leq C^{-j}(-j)!$ for all $j\leq 0$.
}
\eneq
We denote by $\coro$ the subring of $\cor$ consisting of elements of
order $\leq 0$ and by 
$\cor(r)$ the $\coro$-module consisting of elements of order $\leq r$.

We denote by ${}^t(\scbul)\cl\cor\to\cor$
the $\C$-linear automorphism of $\cor$
induced by ${}^t(\tau)=-\tau$ and  call it the transposition. 
We say that a $\C$-linear map $u\cl E\to F$ of 
$\cor$-vector spaces is anti-$\cor$-linear if 
it satisfies $u(a\cdot x)={}^ta\cdot u(x)$ for any 
$x\in E$, $a\in\cor$. 

\subsection*{The ring $\W[T^*X]$.}

On $T^*X$ there exists a {\em no more conic} 
sheaf $\W[T^*X]$ which enjoys the following properties:
\begin{itemize}
\item
$\W[T^*X]$ is a filtered sheaf of central $\cor$-algebras and 
\eqn
&&\gr\W[T^*X]\simeq\sho_{T^*X}[\opb{\tau},\tau].
\eneqn
\item
There is a faithful and flat monomorphism of filtered $\C$-algebras
$\E[T^*X]\hookrightarrow\W[T^*X]$. 
\item
Set 
\eq\label{def:Wv}
&&\Wv[T^*X]=
\opb{\pi}\Omega_X^{\frac{1}{2}}\tens[\opb{\pi}\sho_X]\W[T^*X]
\tens[\opb{\pi}\sho_X]\opb{\pi}\Omega_X^{-\frac{1}{2}}.
\eneq
The sheaf of algebras $\Wv[T^*X]$ is endowed with  an
anti-$\cor$-linear anti-auto\-morphism  $P\mapsto {}^tP$.
\item
Any  symplectic isomorphism
$\psi\cl T^*X\supset U\isoto V\subset T^*Y$
can be {\em locally} quantized as an 
isomorphism of filtered sheaves of $\cor$-algebras commuting with the 
anti-$\cor$-linear anti-isomorphism  $P\mapsto {}^tP$:
\eqn
&& \Psi\cl \oim{\psi}\Wv[T^*X]\isoto \Wv[T^*Y].
\eneqn
(Again, this isomorphism $\Psi$ exists {\em only locally and is not unique}.)
\end{itemize}
Moreover, when $X$ is affine, $\W[T^*X]$ satisfies:
\begin{itemize}
\item
 any section $P\in\W[T^*X](U)$ on an open subset $U\subset T^*X$ admits a 
total symbol
\eq\label{eq:totsymbw}
&&\sigma_{\tot}(P)(x;u,\tau) = \sum_{-\infty<j\leq m}p_j(x;u)\tau^j, 
\,m\in\Z\, \quad p_j\in\sho_{T^*X}(U),
\eneq
with the condition:
\eq\label{eq:defW2}
&&\left\{  \parbox{300 pt}{
for any compact subset $K$ of $U$ there exists a positive 
constant $C_K$ such that 
$\sup\limits_{K}\vert p_{j}\vert \leq C_K^{-j}(-j)!$ for all $j\leq 0$.
}\right. 
\eneq
Note that  $\cor=\W[\rmpt]$.
\item
The total symbol of the product is given by the Leibniz rule:
\eqn
\sigma_{\tot}(P\circ Q)
&=&\sum_{\alpha\in\N^n} \dfrac{\tau^{-\vert\alpha\vert}}{\alpha !} 
\partial^{\alpha}_u\sigma_{\tot}(P)\partial^{\alpha}_x\sigma_{\tot}(Q).
\eneqn 
\item
The total symbol of the transposition is given by 
\eq\label{eq:transposeW}
\sigma_{\tot}({}^tP)(x;u,\tau)
&=&\sum_{\alpha\in\N^n}\dfrac{(-\tau)^{-\vert\alpha\vert}}{\alpha!}
\partial_u^\alpha\partial_x^\alpha{}\sigma_{\tot}(P)(x;u,-\tau).
\eneq
\end{itemize}
We denote by $\W[T^*X](0)$ the subsheaf of 
$\W[T^*X]$ consisting of
sections of order $\leq 0$.
Then  $\W[T^*X](0)$ is a $\coro$-algebra and there is a 
$\cor$-linear isomorphism
\eqn
&&\W[T^*X](0)\tens[\coro_{T^*X}]\cor_{T^*X}\isoto\W[T^*X].
\eneqn

\subsection*{From $\she$ to $\shw$.}

\noindent
One can deduce the algebra $\W[T^*X]$
from the algebra $\E[T^*X]$. 
Let $t\in\C$ be the coordinate and set
$$
\E[T^*(X\times\C),\hat t]=\{P\in\E[T^*(X\times\C)];\ [P,\partial_t] = 0\}.
$$
Set $T^*_{\tau\neq 0}(X\times\C) = \{(x,t;\xi,\tau);\tau\neq 0\}$, and
consider the map
\eqn\label{eq:projrho}
\rho\cl T^*_{\tau\neq 0}(X\times\C) &\to& T^*X,\\
(x,t;\xi,\tau)&\mapsto&(x;\xi/\tau).
\eneqn
The ring $\W[T^*X]$   on $T^*X$ may be defined by setting (see \cite{P-S}):
\eqn
&&\W[T^*X]\eqdot
\oim\rho(\E[T^*(X\times\C),\hat t]\vert_{T_{\tau\neq 0}^*(X\times\C)}).
\eneqn
\begin{remark}
(i) Many authors use the parameter $\hbar$ instead of $\opb{\tau}$. 

\noindent
(ii) There exist formal versions $\HE[T^*X]$ and $\HW[T^*X]$ 
of the sheaves $\E[T^*X]$ and $\W[T^*X]$, respectively,
 and most of the authors work with $\HW[T^*X]$.
\end{remark}

\subsection*
{Deformation quantization and exponential star products.}

\noindent
If $P\in \W[T^*X]$ has order $0$, the operator 
$\exp\tau P$ does not exist in  $\W[T^*X]$.
Using an extra central parameter $t$, 
 a new  
$\cor$-algebra $\TW[T^*X]$ on $T^*X$ was constructed in \cite{Di-S}.
This algebra  enjoys the following properties:

\bnum
\item
there is a monomorphism of
$\cor$-algebras $\iota\cl\W[T^*X]\hookrightarrow \TW[T^*X]$
and a morphism  $\res\cl \TW[T^*X]\to \W[T^*X]$
such that  the composition 
$\W[T^*X]\to \TW[T^*X]\to\W[T^*X]$  is the identity,
\item
any symplectic isomorphism $\psi:T^*X\supset U_X\isoto U_Y\subset T^*Y$
can be locally quantized as an isomorphism of $\cor$-algebras 
$\Psi\cl\TW[T^*X]\isoto\TW[T^*Y]$,
\item
for $P\in\W[T^*X](0)$, the section 
$\exp(t\tau P)=\sum_{n\geq 0}\dfrac{(t\tau P)^n}{n!}$ 
is well defined in $\TW[T^*X]$.
\enum
The algebra $\TW[T^*X]$ is constructed as follows.

Let $s$ be a holomorphic coordinate on $\C$ and 
denote by $\shw_{T^*(\C\times X),\widehat \partial_s}$ 
the subalgebra of $\W[T^*(\C\times X)]$ consisting of sections which
do not depend on $\partial_s$, {\em i.e.,} which commute with $s$. We look
at $\shw_{T^*(\C\times X),\widehat \partial_s}$ as a sheaf on $\C\times T^*X$
and we denote by $p\cl \C\times T^*X\to T^*X$ the projection.
Set 
\eq\label{eq:SW}
&&\SW[T^*X]\eqdot R^1\eim{p}(\shw_{T^*(\C\times X),\widehat \partial_s}).
\eneq
Hence, the sections of $\SW[T^*X]$ are sections of $\W[T^*X]$
depending of an extra holomorphic parameter $s$ defined for 
$\vert s\vert \gg 0$ modulo sections defined for all $s$. 
The convolution product in the $s$ variable allows us to endow 
$\SW[T^*X]$ with a structure of an algebra.

The algebra $\TW[T^*X]$ is constructed as the Laplace transform of 
 $\SW[T^*X]$ which interchanges $s^{-n-1}$ and $\frac{(t\tau)^n}{n!}$.
(Here, $t$ and $\tau$ commute.) 

Note that, if $P$ has order $0$,
the section $s-P$ is invertible on each compact subset $K$ of $T^*X$
 for $\vert s\vert\gg 0$. Therefore, 
$1/(s-P)$ is a well-defined section of $\SW[T^*X]$ and its Laplace
transform $\exp(t\tau P)$ belongs to $\TW[T^*X]$.

\section{Algebroid stacks}\label{section:algebroid}

A local model for a complex symplectic manifold
 is an open subset of $T^*X$. 
Hence, it is natural to ask whether the construction 
of the sheaf of algebras $\shw_{T^*X}$ still makes sense on complex symplectic
manifolds. However, since the quantization of a  symplectic isomorphism is not
unique, one has to replace the notion of
a sheaf of algebras by that of an algebroid stack, a notion introduced in
\cite{Ko}. We refer
to \cite{D-P} for a more systematic study and to \cite{K-S2} 
for an introduction to stacks.

In this section, $\cora$ denotes a commutative unital algebra and $X$ 
a topological space.

If $A$ is a $\cora$-algebra, we denote by $A^+$ the category with 
one object and having $A$ as morphisms of this object. 
Let $\sha$ be a sheaf of $\cora$-algebras on $X$ and consider the prestack 
$U\mapsto \sha(U)^+$ ($U$ open in $X$). 
We denote by $\sha^+$ the associated stack and 
call $\sha^+$ the algebroid stack associated with $\sha$.  

Consider  an open covering $\shu=\{U_i\}_{i\in I}$ of $X$,
sheaves of $\cora$-algebras $\sha_i$ on $U_i$ ($i\in I$) and 
isomorphisms 
$f_{ij}\cl \sha_j|_{U_{ij}} \isoto\sha_i\vert_{U_{ij}}$ ($i,j\in I$).
The existence of a sheaf of $\cora$-algebras $\sha$
locally isomorphic to $\sha_i$ requires the condition $f_{ij}f_{jk} =
f_{ik}$ on triple intersections. 
Let us weaken this last condition by assuming that there exist
invertible sections $a_{ijk} \in \sha_i(U_{ijk})$ satisfying
\eq\label{eq:algebroid}
&& \left\{ \parbox{300 pt}{
$f_{ij}f_{jk} = \rmad(a_{ijk}) f_{ik}$ on $U_{ijk}$, \\
$a_{ijk} a_{ikl} = f_{ij}(a_{jkl}) a_{ijl}$ on $U_{ijkl}$.
}\right. 
\eneq
(Recall that $\rmad(a)(b)=a\cdot b\cdot\opb{a}$.)
One calls 
\eq\label{eq:descentdata}
&&( \{\sha_i\}_{i\in I},\{f_{ij}\}_{i,j\in I}, \{a_{ijk}\}_{i,j,k \in I})
\eneq
a  descent datum for $\cora$-algebroid stacks on $\shu$.
For such a descent datum, we shall denote by 
$f^+_{ij}\cl \sha_j^+\isoto\sha_i^+$
the equivalences of stacks associated with the isomorphisms $f_{ij}$. 
The following result is stated (in a different form) in \cite{Ka1} and 
goes back to \cite{Gi}.
\begin{theorem}\label{th:descdata}
Consider a 
descent datum \eqref{eq:descentdata} on $\shu$. Then there exist 
a stack $\sha^+$ on $X$, equivalences of stacks 
$\phi_i\cl \sha^+\vert_{U_i}\isoto \sha_i^+$ and isomorphisms of functors
 $c_{ij}\cl f^+_{ij}\isoto \phi_{i}\circ\opb{\phi_j}$
satisfying $c_{ij} \circ c_{jk} \circ a_{ijk} = c_{ik}$.
Moreover, the data $(\sha^+,\{\phi_i\}_i,\{c_{ij}\}_{ij})$ are unique
up to equivalence of stacks, this equivalence being unique up to
a unique isomorphism.
\end{theorem}
One calls $\sha^+$ an algebroid stack.
Although $\sha^+$ is not a sheaf of algebras, modules over $\sha^+$ are
well defined. They are described by pairs
\eqn
&&\shm = ( \{\shm_i\}_{i \in I}, \{\xi_{ij}\}_{i,j \in I} ),
\eneqn
where $\shm_i$ is an $\sha_i$-module and
$\xi_{ij}\cl \atw{f_{ji}}\shm_j\vert_{U_{ij}} \to \shm_{i}
\vert_{U_{ij}}$ is an isomorphism of $\sha_i$-modules
 such that for any $u_k\in\shm_k$ one has
\begin{equation}
\label{eq:glueM}
\xi_{ij}  ({}_{f_{ji}}\xi_{jk} (u_k))= \xi_{ik}(a_{kji}^{-1}  u_k).
\end{equation}
Here, $\atw{f_{ji}}\shm_j$ is the $\sha_i$-module deduced from the 
$\sha_j$-module $\shm_j\vert_{U_{ij}}$ by the isomorphism $f_{ji}$.

One gets a Grothendieck category $\md[\sha^+]$ and the
prestack $\Mods(\sha^+)$ given by 
$U\mapsto \md[\sha^+\vert_U]$ is a stack equivalent
on $U_i$ to the stack $\Mods(\sha_i)$.

\subsection*{Twisted sheaves.}

As a particular case of a module over an algebroid stack, one has the
notion of a twisted sheaf. Assume that $\cora$ is a field and denote
by $\cora^\times$ the group of its invertible elements.

Let $X$ be  a manifold and let 
${\bold c}\in H^2(X;\cora^\times)$.
Represent ${\bold c}$ by a \v{C}ech cocycle 
$\{c_{ijk}\}_{i,j,k \in I}$ associated to an open covering 
$\shu=\{U_i\}_{i\in I}$ of $X$.
We thus get a  descent datum for $\cora$-algebroid stacks
\eqn
&&\cora_{X,{\bold c}}\eqdot ( \{ \cora_{U_i} \}_{i \in I}, 
\{\id_{\cora_{U_{ij}} } \}_{i,j \in I}, \{c_{ijk}\}_{i,j,k \in I}),
\eneqn
and cohomologous cocycles give equivalent stacks.

\begin{example}\label{ex:Omega1/2}
Assume now that $X$ is a complex manifold. 
Consider the short exact sequence
\eqn
&&1 \to \C_X^\times \to \OO_X^\times \to[d\log] d\OO_X \to 0
\eneqn
which gives rise to the 
long exact sequence
\eqn
&&H^1(X;\C_X^\times) \to[\alpha]H^1(X;\OO_X^\times) \to[\beta]
H^1(X;d\OO_X) \to[\gamma]H^2(X;\C_X^\times).
\eneqn
If $\shl$ is  a line bundle, it defines a class 
$[\shl]\in H^1(X;\OO_X^\times)$. For $\lambda\in\C$, one sets
\eqn
&& {\bf c}_{\shl}^\lambda=
\gamma(\lambda\cdot\beta([\shl]))\in H^2(X;\C_X^\times).
\eneqn
We shall apply this construction when $\shl=\Omega_X$ and 
$\lambda=\frac12$ and set for short:
\eqn
&& \md[\C_{X,\frac12}]=\md[\C_{X,{\bf c}_{\Omega_X}^\frac12}].
\eneqn
\end{example}

\section{Quantization of symplectic manifolds}

On any complex contact manifold, 
the existence of a canonical $\C$-algebroid stack locally equivalent to the algebroid 
stack associated with the sheaf of algebras of microdifferential operators of \cite{S-K-K}
has been obtained by M. Kashiwara in \cite{Ka1}.

On any complex Poisson manifold, the existence of 
a $\hcor$-algebroid stack of formal deformation quantization
has been obtained by M. Kontsevich \cite{Ko}. The analytic case
on symplectic manifolds has been obtained in \cite{P-S} by
a different method, making a link with Kashiwara's construction.
The classification of these algebroid stacks is discussed in \cite{Po}.

In particular, for a complex symplectic manifold $\stx$, 
there is a canonical $\cor$-algebroid stack $\Wva[\stx]$ locally equivalent to
the algebroid stack $\Wva[T^*X]$ associated with the sheaf of
algebras $\Wv[T^*X]$. 
The same result holds with $\W[\stx]$ replaced by $\W[\stx](0)$ 
and $\cor$ by $\coro$. 
\begin{notation}
For short, as far as there is no risk of confusion, 
we shall write $\W[\stx]$ instead of $\Wva[\stx]$.
\end{notation}

Let $\stx$ be a complex symplectic manifold.
Then $\md[{\W[\stx]}]$ is a Grothendieck category.
We denote by $\RD^\Rb(\W[\stx])$ its bounded derived category and
call an object of this derived category a
$\W[\stx]$-module. 
One proves as usual that the sheaf of algebras $\W[T^*X]$ is
coherent and  the support
of a coherent $\W[T^*X]$-module is a closed complex analytic
subvariety of $T^*X$. This support is involutive in view of 
Gabber's theorem (see \cite[Th.~7.33]{Ka2}). Hence, the (local) notions of a 
coherent or holonomic $\W[\stx]$-module make sense. 

Similarly as for $\shd$-modules (see \cite{Ka2}), one says that a 
coherent $\W[\stx]$-module $\shm$ is good if, for any
open relatively compact subset $U$ of $\stx$, there exists 
a coherent $\W[\stx](0)\vert_U$-module $\shm_0$ contained in  $\shm\vert_U$ 
which generates $\shm\vert_U$.

Let us denote by:
\begin{itemize}
\item
$\RD^\Rb_{\coh}(\W[\stx])$ the full triangulated subcategory 
 of $\RD^\Rb(\W[\stx])$ consisting
of objects with coherent cohomologies.
\item
$\RD^\Rb_{\gd}(\W[\stx])$ the full triangulated subcategory of 
$\RD^\Rb_{\coh}(\W[\stx])$ consisting of objects with good cohomologies. 
\item
$\RD^\Rb_{\gdc}(\W[\stx])$ the full triangulated subcategory of 
$\RD^\Rb_{\gd}(\W[\stx])$ consisting of objects with compact supports.
\item 
$\RD^\Rb_{\hol}(\W[\stx])$ the full triangulated subcategory 
 of $\RD^\Rb_{\coh}(\W[\stx])$ consisting 
of objects with
Lagrangian supports in $\stx$. (One calls such an object an holonomic 
$\W[\stx]$-module.)
\item 
$\RD^\Rb_{\rh}(\W[\stx])$ the full triangulated 
subcategory of $\RD^\Rb_{\hol}(\W[\stx])$ consisting of objects with
regular holonomic cohomologies (to be defined below).
\item
Let $\stx$ and $\sty$ be two complex symplectic manifolds and 
let $\shm\in\RD^\Rb(\W[\stx])$, 
$\shn\in\RD^\Rb(\W[\sty])$. Their exterior product is given by
$\shm\detens\shn\eqdot\W[\stx\times\sty]
\etens_{\W[\stx]\etens\W[\sty]}(\shm\etens\shn)$. 
\end{itemize}

\subsection*{Simple $\W[\stx]$-modules.}

\begin{definition}\label{def:simpleWmod}
Let $\Lambda$ be a smooth Lagrangian submanifold of $\stx$. 
\banum
\item
Let $\shl(0)$ be a coherent $\W[\stx](0)$-module supported by $\Lambda$. 
One says that $\shl(0)$ is simple along $\Lambda$ if
$\shl(0)/\shl(-1)$ is an invertible $\OO_\Lambda$-module. 
Here,  $\shl(-1)=\cor_\stx(-1)\shl(0)$.
\item
Let $\shl$ be a coherent $\W[\stx]$-module supported by $\Lambda$.  
One says that $\shl$ is simple along
$\Lambda$ if there locally exists a coherent $\W[\stx](0)$-submodule $\shl(0)$
of $\shl$ such that $\shl(0)$  generates $\shl$ over $\W[\stx]$ and is 
simple along $\Lambda$.
\item
Let $\shl$ be a coherent $\W[\stx]$-module supported by $\Lambda$.  
One says that $\shl$ is regular if, locally, it is a finite direct 
sum of simple modules.
\item 
Let $\Lambda$ be a, not necessarily 
smooth, Lagrangian subvariety of $\stx$. A
coherent $\W[\stx]$-module supported by $\Lambda$
 is regular if it is regular at generic points of $\Lambda$.
One calls such an object a regular holonomic $\W[\stx]$-module.
\eanum
\end{definition}
It follows from Gabber's theorem that when $\Lambda$ is smooth, 
Definitions~\ref{def:simpleWmod} (c) and (d) coincide 
(see \cite[Th. 8.34]{Ka2}).

One proves easily that 
any two $\W[\stx]$-modules simple along $\Lambda$ are locally
isomorphic and that if $\shl_i$ ($i=0,1$) are simple along 
$\Lambda$, 
then $\rhom[{\W[\stx]}](\shl_0, \shl_1)$ is concentrated in degree $0$ and is 
a $\cor$-local system of rank one on $\Lambda$.

\begin{example}
Let $X$ be  a complex manifold. 
We denote by $\OO^\tau_X$ the 
$\W[T^*X]$-module supported by the zero-section $T^*_XX$ defined 
by $\OO^\tau_X = \W[T^*X]/\shi$, where $\shi$ is the left ideal 
generated by the vector fields which annihilate the section
$1\in\sho_X$. A section $f(x,\tau)$
of this module may be written as a series:
\eq
&& f(x,\tau)=\sum_{-\infty<j\leq m}f_j(x)\tau^j,\quad m\in\Z,
\eneq
the $f_j$'s satisfying Condition \eqref{eq:defW2}.
Then $\OO^\tau_X$ is a simple $\W[T^*X]$-module along $T^*_XX$. 
\end{example}
 
The next result asserts that, up to a twist, there exist globally
defined simple $\W[\stx]$-modules. 
 
\begin{theorem} {\rm \cite{Da-S}}
Let $\Lambda$ be a smooth Lagrangian submanifold of $\stx$.
There is an equivalence of $\cor$-additive stacks:
\eq\label{eq:DSeqv}
\Mods_{\text{reg-}\Lambda}(\W[\stx])|_\Lambda \simeq
\Mods_{\text{loc-sys}}(\cor_\Lambda\tens[\C]\C_{\Lambda,1/2}).
\eneq
\end{theorem}
Here, the  left-hand side is the substack of 
$\Mods(\W[\stx])|_\Lambda$ 
consisting of regular holonomic modules  along $\Lambda$ 
and the right-hand side is the substack of the stack of 
twisted sheaves of $\cor_\Lambda$-modules with twist 
$\C_{\Lambda,1/2}$ consisting of objects locally isomorphic 
to local systems over $\cor$.
The proof uses the corresponding theorem
 for contact manifolds  due to Kashiwara.

\subsection*{$\W[\stx]$-module associated with the diagonal.}

Let $\stx$ be a complex symplectic manifold. We  
denote
by $\stx^a$ the complex manifold $\stx$ endowed with 
the symplectic form  $-\omega$, where  $\omega$ is the symplectic form 
on $\stx$.
There is a natural equivalence of
algebroid stacks $\W[\stx]^\rop\simeq\W[\stx^a]$. 

We denote  by $\Delta_\stx$ the diagonal of $\stx\times\stx^a$ and by  
$d_\stx$ the complex dimension of $\stx$.
\begin{theorem}\label{th:existCD}
\bnum
\item
There exists a simple 
 $\W[\stx\times\stx^a]$-module $\CD[\Delta_\stx]$
supported by the diagonal $\Delta_\stx$ of  $\stx\times\stx^a$ with
the property that if $U$ is open in $\stx$ and isomorphic to an open
subset $V$ of a cotangent bundle $T^*X$, then $\CD[\Delta_\stx]\vert_U$
is isomorphic to $\W[T^*X]\vert_V$ as a $\W[T^*X]\tens\W[T^*X]^\rop$-module. 
\item
There is a  natural isomorphism
$\CD[\Delta_{\stx^a}]\lltens[{\W[\stx\times\stx^a]}]\CD[\Delta_\stx]
\simeq \cor_{\Delta_\stx}\,[d_\stx]$.
\enum
\end{theorem}
(i) follows from general considerations on algebroid stacks. (ii)
follows from a construction of Feigin and Tsygan \cite{F-T} 
(see also \cite{E-F}).

Let $\shm\in\RD^\Rb(\W[\stx])$. We set
\eq
&& \RD'_{\rmw}\shm\eqdot \rhom[{\W[\stx]}](\shm,\CD[\Delta_\stx]),
\quad \RD_{\rmw}\shm\eqdot \RD'_{\rmw}\shm\,[\frac{1}{2}d_\stx].
\eneq
These objects are well-defined in $\RD^\Rb(\W[\stx^a])$.
Using Theorem~\ref{th:existCD} and the fact that $\CD[\Delta_\stx]$
is simple, 
one gets  an isomorphism in $\RD^\Rb_{\gd}(\W[\stx^a\times\stx])$:
\eq\label{eq:HHisom}
&&\RD_\rmw(\CD[\Delta_\stx])\simeq\CD[\Delta_{\stx^a}].
\eneq
\begin{remark}\label{rem:rigdual}
By extending the definition of Van den Bergh \cite{VdB2} 
(see also \cite{Y-Z}) to algebroid
stacks, the isomorphism \eqref{eq:HHisom} could be translated by saying that 
$\CD[\Delta_\stx]\,[d_\stx]$ is a rigid dualizing complex over $\W[\stx]$.
\end{remark}
Let $\shm,\shn$ be two objects of 
$\RD^\Rb_{\coh}(\W[\stx])$. Then, after identifying $\Delta_\stx$
with $\stx$ by the first projection, there is a natural isomorphism
in $\RD^\Rb(\cor_\stx)$:
\eqn
&&\rhom[{\W[\stx]}](\shm,\shn)\simeq 
\rhom[{\W[\stx\times\stx^a]}]
(\shm\detens\RD'_{\rmw}\shn,\CD[\Delta_\stx]).
\eneqn

\section{Constructibility and perversity}

We refer to \cite{K-S1} for basic notions on sheaves. 
Let $Z$ be  a real analytic manifold. Recall that one denotes by:
\begin{itemize}
\item
$C(S_1,S_2)$ the normal cone of two subsets $S_1$ and $S_2$ of $Z$, a
closed conic subset of the tangent space $TZ$, identified with a
subset of $T^*Z$ in case $Z$ is symplectic,
\item
 $\RD^\Rb(\cor_Z)$ the bounded derived category of
sheaves of $\cor$-modules on $Z$,
\item
$\SSi(F)$ the microsupport of an object $F\in \RD^\Rb(\cor_Z)$, a
closed conic involutive subset of $T^*Z$,
\item
 $\ori_Z$ the orientation sheaf on $Z$, $\omega_Z$
the dualizing complex (hence, $\omega_Z\simeq \ori_Z\,[\dim_\R Z]$),
$\RD_Z\eqdot \rhom[\cor_Z](\scbul,\omega_Z)$
 the duality functor for sheaves,
\item
$\RD_\Rc^\Rb(\cor_Z)$
the full triangulated subcategory of $\RD^\Rb(\cor_Z)$ consisting of
objects with $\R$-constructible cohomology and, in case $Z$ is a
complex manifold, $\RD_\Cc^\Rb(\cor_Z)$ 
the full triangulated subcategory of $\RD^\Rb(\cor_Z)$ consisting of
objects with $\C$-constructible cohomology.
\end{itemize}

\begin{theorem}\label{th:main1}{\rm \cite{K-S3}}
Let $\stx$ be a complex symplectic manifold and let 
$\shl_i$ \lp$i=0,1$\rp\, be two objects of 
$\RD^\Rb_{\rh}(\W[\stx])$ supported by smooth Lagrangian 
manifolds  $\Lambda_i$. Then
\bnum
\item
the object $\rhom[{\W[\stx]}](\shl_1,\shl_0)$ 
belongs to  $\RD^\Rb_\Cc(\cor_\stx)$ and its microsupport 
is contained in the normal cone $C(\Lambda_0,\Lambda_1)$,
\item
the natural morphism
\eqn
&&\rhom[{\W[\stx]}](\shl_1,\shl_0)\to
\RD_\stx\bigl(\rhom[{\W[\stx]}](\shl_0,\shl_1\,[d_\stx])\bigr)
\eneqn
is an isomorphism.
\enum
\end{theorem}
The proof makes use of tools from the theory of holonomic 
$\shd$-modules and uses
some functional analysis, namely Houzel's theorem \cite[2,\S~4~Th.~1']{Ho}.
\begin{remark}\label{rem:kai}
In \cite{B-F}, K. Behrend and B. Fantecci
construct complexes (over $\C$) 
naturally associated with the data of two smooth Lagrangian
submanifolds. Their result should have some relations with 
Theorem~\ref{th:main1}.
 \end{remark}
\begin{corollary}\label{cor:main1}
Let  $\shl_0$ and $\shl_1$ be two regular holonomic 
$\W[\stx]$-modules supported by smooth Lagrangian manifolds.
Then the object 
$\rhom[{\W[\stx]}](\shl_1,\shl_0)$ of $\RD^\Rb_\Cc(\cor_\stx)$
is perverse.
\end{corollary}

\begin{conjecture}{\rm \cite{K-S3}}\label{conj:mainconj}
Theorem~\ref{th:main1} remains true without assuming that 
the $\Lambda_i$'s are smooth.
\end{conjecture}

\section{Composition of kernels and Calabi-Yau categories}

Consider three complex symplectic manifolds $\stx_i$ ($i=1,2,3$) and 
denote as usual by $p_i$ and 
$p_{ji}$ the projections defined on $\stx_3\times\stx_2\times\stx_1$. 

For $\Lambda_i$ a closed subset of $\stx_{i+1}\times\stx_i$ ($i=1,2$), we set
\eq\label{eq:compLambda}
&&\Lambda_2\circ\Lambda_1\eqdot 
p_{31}(\opb{p_{32}}\Lambda_2\cap\opb{p_{21}}\Lambda_1).
\eneq
For $\shk_i\in\RD^\Rb_{\gd}(\W[\stx_{i+1}\times\stx_i^a])$ ($1\leq i\leq 2$),
we set
\eq\label{eq:compKer}
&& \shk_2\circ\shk_1\eqdot 
\reim{p_{31}}
(\opb{p_{32}}\shk_2\lltens[\opb{p_2}{\W[\stx_2]}]\opb{p_{21}}\shk_1).
\eneq

\begin{theorem}\label{th:compkern}{\rm \cite{Sc-S2}}
Assume that $p_{31}$ is proper on 
$\opb{p_{32}}\supp(\shk_2)\cap\opb{p_{12}}\supp(\shk_1)$. 
Then 
\bnum
\item
the object $\shk_2\circ\shk_1$ belongs to 
$\RD^\Rb_{\gd}(\W[\stx_3\times\stx_1^a])$,
\item
there is a natural  isomorphism
in $\RD^\Rb_{\gd}(\W[\stx_3^a\times\stx_1])$:
\eq\label{eq:dualitymor1}
&&\RD_{\rmw}\shk_2\circ\RD_{\rmw}\shk_1\isoto\RD_{\rmw}(\shk_2\circ\shk_1).
\eneq
\enum
\end{theorem}
The proof of~(i) uses again \cite[2,\S~4~Th.~1']{Ho}.
The construction of the duality morphism in~(ii) 
uses the isomorphism \eqref{eq:HHisom}.

Choosing $\stx_3=\stx_1=\rmptt$ in Theorem~\ref{th:compkern},
we get:

\begin{corollary}\label{co:compkern1}
Let $\stx$ be a complex symplectic manifold and let 
$\shm$ and $\shn$ be two objects of 
$\RD^\Rb_{\gd}(\W[\stx])$. 
Assume that $\supp(\shm)\cap\supp(\shn)$ is compact. Then
\bnum
\item
the object $\RHom[{\W[\stx]}](\shm,\shn)$ 
 has $\cor$-finite dimensional cohomology,
\item
there is a natural isomorphism, functorial with respect to $\shm$ and $\shn$:
\eqn
&&\RHom[{\W[\stx]}](\shm,\shn)\simeq
\bigl(\RHom[{\W[\stx]}](\shn,\shm\,[d_\stx])\bigr)^\star,
\eneqn
where ${}^\star$ is the duality functor for $\cor$-vector spaces.
\enum
\end{corollary}

Recall \cite{B-K} that a  $\cor$-triangulated category $\sht$
is $Ext$-finite if for any two objects $F,G$ of $\sht$, 
the $\cor$-vector space $\bigoplus_{i\in\Z}\Hom[\sht](F,G[i])$ 
is finite dimensional.
In this situation, a Serre functor $S\cl\sht\to\sht$ is an equivalence
of $\cor$-triangulated categories such that
\eqn
&&(\Hom[\sht](F,G))^\star\simeq \Hom[\sht](G,S(F))
\eneqn
functorially in $F$ and $G$.  If the Serre functor is a shift by an 
integer $d$, one
says that $\sht$ is a Calabi-Yau category of dimension $d$.

\begin{corollary}\label{co:compkern2}
Let $\stx$ be a complex symplectic manifold. 
Then 
$\RD^\Rb_{\gdc}(\W[\stx])$ is a Calabi-Yau category of dimension $d_\stx$.
\end{corollary}

\begin{remark}
(i) The analogue of Corollary~\ref{co:compkern1}  
on complex contact manifolds over the field $\C$ is false.
Note that Corollary~\ref{co:compkern1}  may be considered 
as a direct image theorem over one point, 
and a manifold of dimension $0$  has a complex  symplectic structure,
not a complex contact structure.

Nevertheless, the analogue of Corollary~\ref{co:compkern1} 
is true over the field $\C$ for
complex contact manifolds when restricting to the category of
regular holonomic modules. This follows from 
results obtained in \cite{K-K}.

\noindent
(ii) On a complex compact manifold $X$, 
the bounded derived category of sheaves with $\C$-constructible cohomology
is an $Ext$-finite triangulated category,
as well as the equivalent category, the bounded derived 
category of $\shd_X$-modules  with regular holonomic cohomology.
Both categories do not seem to have a Serre functor.
\end{remark}

\section{Index theorem}
In this section, we announce works in progress with
J-P.~Schneiders~\cite{Sc-S2}. 

\subsection*{Euler class.}
Let $\stx$ be complex symplectic manifold and let
$\shm\in\RD^\Rb_{\coh}(\W[\stx])$. 
We have the chain of morphisms
\eqn
\rhom[{\W[\stx]}](\shm,\shm)
&\isofrom&\rhom[{\W[\stx]}](\shm,\CD[\Delta_\stx])
\lltens[{\W[\stx]}]\shm\\
&\simeq&
(\rhom[{\W[\stx]}](\shm,\CD[\Delta_\stx])\tens[\cor_{\stx}]\shm)
\lltens[{\W[\stx]\tens\W[\stx^a]}]\CD[\Delta_\stx]\\
&\to&\CD[\Delta_{\stx^a}]\lltens[{\W[\stx\times\stx^a]}]\CD[\Delta_\stx]
\simeq\cor_\stx\,[d_\stx].
\eneqn
Here, we have used Theorem~\ref{th:existCD}~(ii).
We get a map
\eqn
&&\Hom[{\W[\stx]}](\shm,\shm)\to H^{d_\stx}_{\supp(\shm)}(\stx;\cor_\stx).
\eneqn
The image of $\id_\shm$ gives an element 
\eq\label{eq:defeuM}
&&\Eu(\shm)\in H^{d_\stx}_{\supp(\shm)}(\stx;\cor_\stx).
\eneq
\subsection*{Symplectic Riemann-Roch theorem}

We consider the situation of Theorem~\ref{th:compkern}. Hence, we have
three complex symplectic manifolds $\stx_i$ ($i=1,2,3$) and we have 
closed subsets $\Lambda_i$ of $\stx_{i+1}\times\stx_i$ ($i=1,2$).
We set for short $d_i\eqdot d_{\stx_i}$ ($i=1,2,3$) and consider
cohomology classes 
$\lambda_i\in H^{d_{i+1}+d_i}_{\Lambda_{i+1}\times\Lambda_i}(\stx_{i+1}\times\stx_{i};
\cor_{\stx_{i+1}\times\stx_{i}})$ ($i=1,2$).
Assuming that $p_{31}$ is proper on 
$\opb{p_{32}}(\Lambda_2)\cap\opb{p_{12}}(\Lambda_1)$,
 we set
\eqn
\lambda_2\circ\lambda_1&\eqdot& 
\int_{\stx_2}(\opb{p_{32}}\lambda_2\cup\opb{p_{21}}\lambda_1)\in 
H^{d_3+d_1}_{\Lambda_2\circ\Lambda_1}(\stx_3\times\stx_1;
\cor_{\stx_3\times\stx_1 }).
\eneqn
Here, $\cup$ is the cup product and 
\eqn
&&\int_{\stx_2}\cl 
H^{d_3+2d_2+d_1}_{\opb{p_{32}}\Lambda_2\cap\opb{p_{21}}\Lambda_1}
(\stx_3\times\stx_2\times\stx_1;\cor_{\stx_3\times\stx_2\times\stx_1 })\to
H^{d_3+d_1}_{\Lambda_2\circ\Lambda_1}(\stx_3\times\stx_1;
\cor_{\stx_3\times\stx_1 })
\eneqn
is the Poincar{\'e} integration morphism. 

Let $\shk_i\in\RD^\Rb_{\gd}(\W[\stx_{i+1}\times\stx_i^a])$ \lp$1\leq i\leq 2$\rp\,
and assume that $p_{31}$ is proper on 
$\opb{p_{32}}\supp(\shk_2)\cap\opb{p_{12}}\supp(\shk_1)$. 
The proof of the formula 
\eq\label{eq:euRR1}
&&\Eu(\shk_2\circ\shk_1)=
\Eu(\shk_2)\circ \Eu(\shk_1)
\eneq
is in progress. It would partly generalize the index theorems  for
coherent  $\shd_X$-modules proved in \cite{Sc-S}.

As a particular case of \eqref{eq:euRR1}, one finds that for two objects 
$\shl$ and $\shm$ in $\RD^\Rb_{\gd}(\W[\stx])$ such that 
$\supp\shl\cap\supp\shm$ is compact, we have
\eq\label{conj:RR1}
&&\chi(\RHom[{\W[\stx]}](\shl,\shm))
=\int_{\stx}\Eu(\RD'_{\rmw}\shl)\cup \Eu(\shm).
\eneq
In the case of coherent  $\shd_X$-modules on a complex manifold $X$,  
the formula
\eqn
&&\Eu(\shl)=[\Ch(\gr\shl)\cup \Td(TX)]^{d_{T^*X}}
\eneqn
had been conjectured in  \cite{Sc-S} and proved by \cite{B-N-T}. 
On a complex symplectic manifold, these authors
 give a more general formula  calculating $\Eu(\shl)$, a formula
in which the class of the deformation quantization appears.

\end{document}